\documentclass[12pt]{article}

\usepackage{amssymb}

\usepackage[utf8]{inputenc}
\usepackage[T1]{fontenc}
\usepackage{pgf,pgfarrows,pgfnodes,pgfautomata,pgfheaps}
\usepackage{amscd}
\usepackage{amsfonts}
\usepackage{latexsym}
\usepackage{graphicx}
\usepackage{amsthm}
\usepackage{amsmath}
\usepackage{leftidx}
\usepackage{wasysym}
\usepackage{wrapfig}

\def\Z{\mathbb{Z}}

\def\Q{\mathbb{Q}}
\def\C{\mathbb{C}}
\def\F{\mathbb{F}}

\def\det{\mathrm{det}}

\def \det{\operatorname{det}}

\theoremstyle{definition}

\newtheorem{thm}{Theorem}
\newtheorem{df}[thm]{Definition}
\newtheorem{cor}[thm]{Corollary}

\newtheorem{lem}[thm]{Lemma}
\newtheorem{rmk}[thm]{Remark}
\newtheorem{ex}[thm]{Example}

\newcommand\myeqq{\mathrel{\overset{\makebox[0pt]{\mbox{\normalfont\tiny\sffamily lemma 6}}}{=}}}

\begin{document}

\date{}
\author{}

 \let\thefootnote\relax\footnote{ This work was partially supported by the Faculty of Applied Mathematics AGH UST statutory tasks within subsidy of Ministry of Science and Higher Education.
 
 This work was partially supported by the Faculty of Mathematics and Computer Science at Jagiellonian University statutory tasks within subsidy of Ministry of Science and Higher Education.
 
 This research was supported in part by PLGrid Infrastructure.}

\begin{center}

{\Large \textbf{Algorithm for studying polynomial maps\\ and reductions modulo prime number}}

\vspace{20pt}

EL\.ZBIETA ADAMUS \\ Faculty of Applied Mathematics, \\ AGH University of Science and Technology \\
al. Mickiewicza 30, 30-059 Krak\'ow, Poland \\
e-mail: esowa@agh.edu.pl \\

\vspace{0.5cm}
PAWE\L \ BOGDAN \\ Faculty of Mathematics and Computer Science, \\ Jagiellonian University \\
ul. \L ojasiewicza 6, 30-348 Krak\'ow, Poland \\
e-mail: pawel.bogdan@uj.edu.pl \\

\end{center}
\vspace{10pt}

\begin{abstract}
In \cite{ABCH} an effective algorithm for inverting polynomial automorphisms was proposed. Also the class of Pascal finite polynomial
automorphisms was introduced. Pascal finite polynomial maps constitute a generalization of exponential automorphisms to positive characteristic. 
 In this note we explore properties of  the algorithm while using Segre homotopy and reductions modulo prime number. We give a method of retrieving an inverse of a given polynomial automorphism $F$ with integer coefficients form a finite set of the inverses of its reductions modulo prime numbers.
Some examples illustrate effective
aspects of our approach.

\end{abstract}

\section{Introduction}

Let $K$ be a field and
let $F=(F_1, \ldots, F_n) :K^n \rightarrow K^n$ be a
polynomial map. $F$ is invertible over $K$ if there exists a polynomial mapping $G : K^n \rightarrow K^n$ such that $F \circ  G = Id$ and $G \circ F = Id$. 
Study of invertible polynomial mappings is related to the famous Jacobian Conjecture, which asks if every
polynomial mapping such that its jacobian is nonzero constant is invertible with polynomial inverse. Many results concerning polynomial automorphisms are formulated for an arbitrary field $K$, but the case of a field of characteristic zero is the one discussed most often. However after reducing coefficients of $F \in \mathbb{Z}[X]^n$ modulo given prime number one can consider it over finite field $\mathbb{F}_p$. Results concerning this approach can be found for example in \cite{MFF}, \cite{MW}.

In \cite{ABCH} we described an algorithm which for a given $F \in K[X]^n$ over an arbitrary field $K$ constructs recursively a sequence of polynomial maps.
We define an endomorphism $\sigma_F$ of $K[X]^n$ 
by $\sigma_F(P)=P\circ F$ and a $\sigma_F$-derivation $\Delta_F$ on $K[X]^n$ by $\Delta_F(P)=\sigma_F(P)-P$.
Following Maple environment commands we use term \emph{lower degree} instead of an order of vanishing of a polynomial.
We consider $F\in K[X]^n$ of the form
  \begin{equation}\label{xh}
      \left\{ \begin{array}{lll}
             F_1(X_1, \ldots, X_n) &=& X_1+H_1(X_1, \ldots, X_n)\\
              & \vdots & \\
             F_n(X_1, \ldots, X_n) &=& X_n+H_n(X_1, \ldots, X_n),
            \end{array} \right. 
  \end{equation}
  where $H_i(X_1, \ldots, X_n)$ is a  polynomial in $X_1, \ldots, X_n$ of degree $D_i$ and lower degree $d_i$, with $ d_i\geq 2$, for $i=1, \ldots, n$. Let $d=\min d_i, D=\max D_i$.
Then we consider the sequence $P_l=(P_l^1, \ldots, P_l^n)$ of
polynomial maps in $K[X]^n$ defined by $P_{l}=\Delta^l_F(Id)$, where $Id(X)=X$ and $\Delta^l_F$ denotes $\Delta_F \circ \stackrel{l}{\dots} \circ \Delta_F$.
The class of polynomial automorphisms for which the algorithm stops has been distinguished. 
Polynomial map $F:K^n \rightarrow K^n$ is called \emph{Pascal finite} if there exists $m$ such that $P_m=0$. Then $F$ is invertible and the inverse map $G$ of $F$ is given by
 \begin{equation}\label{inveq}
     G(X)= \sum_{l=0}^{m-1} (-1)^l P_l(X)
 \end{equation}
 (see \cite{ABCH}, corollary 2.1).
Pascal finite automorphisms are roots of a polynomial of the form $P(X) = (X-1)^m$. In \cite{ABCH2} we discussed their properties.
They are natural generalization of exponential automorphisms
to positive characteristic.

 In this paper we consider polynomial maps over $\mathbb{Q}$. Those can be transformed into maps with coefficients in $\mathbb{Z}$ by using Segre homotopy also known as denominators clearing procedure. 
 Using clearing map Connel and van den Dries proved (see \cite{CD}, theorem 1.5 or \cite{E}, proposition 1.1.19) that if there is a counterexample to the Jacobian Conjecture $h:\C^m\rightarrow \C^m$, then for some $n >m$ there is a counterexample $f: \C^n \rightarrow \C^n$ with coefficients in $\Z$. In fact Jacobian Conjecture over $\mathbb{C}$ is equivalent to the Jacobian Conjecture over $\mathbb{Z}$ (see \cite{E}, proposition 1.1.12). That is why one can be interested in studying maps with integer coefficients.
 We discuss behaviour of the algorithm while using Segre homotopy. After that we perform reduction modulo prime number $p$ and apply the algorithm proposed in \cite{ABCH}  in order to find an inverse of a reduced map.
 We explore a method of retrieving an inverse of a given polynomial automorphism $F\in \mathbb{Z}[X]^n$ from a finite set of the inverses of its reductions modulo prime numbers. 

Below we recall the main result of \cite{ABCH} (see theorem 3.1) which formulates an equivalent condition to invertibility of a polynomial map and explains how Pascal finite automorphisms admit an algorithmic treatment. This theorem allows to to check if a given polynomial map is invertible
and to find an inverse of a given polynomial automorphism even if it is not a Pascal finite one.

 \begin{thm}\label{symthm}

   Let $F=(F_1, \ldots, F_n) :K^n \rightarrow K^n$ be a
polynomial map of the form (\ref{xh}).
  The following conditions are equivalent:

\begin{enumerate}
\item $F$ is invertible.
\item For $i=1,\ldots, n$ and every $m > \frac{D^{n-1}-d_i}{d-1}+1$, we have
   \begin{equation}
    \sum_{j=0}^{m-1}(-1)^j P^i_j(X)=G_i(X)+R^i_m(X),
    \label{symmetry}
   \end{equation}
 
   where $G_i(X)$ is a polynomial of degree $ \leq D^{n-1}$, independent of $m$, and $R^i_m(X)$ is a polynomial satisfying $R^i_m(F)=(-1)^{m+1}P^i_m(X)$
   (with lower degree $\geq (m-1)(d-1)+d_i>D^{n-1}$).
\item For $i=1,\ldots, n$ and  $m = \lfloor \frac{D^{n-1}-d_i}{d-1}+1 \rfloor +1$, we have
   \[ \sum_{j=0}^{m-1}(-1)^j P^i_j(X)=G_i(X)+R^i_m(X). \]
   where $G_i(X)$ is a polynomial of degree $ \leq D^{n-1}$,
  and $R^i_m(X)$ is a polynomial satisfying $R^i_m(F)=(-1)^{m+1}P^i_m(X)$.
\end{enumerate}

  Moreover the inverse $G$ of $F$ is given by
\begin{equation}
    G_i(Y_1, \ldots, Y_n)=\sum_{l=0}^{m-1}(-1)^l\widetilde{P}^i_l(Y_1, \ldots, Y_n), \, i=1, \ldots, n,
    \label{thminveq}
\end{equation}
  where $\widetilde{P}^i_l$ is the sum of homogeneous summands of $P^i_l$ of degree $ \leq D^{n-1}$ and $m$ is an integer $> \frac{D^{n-1}-d_i}{d-1}+1$.

\label{thm11}
 \end{thm}

 \section{Segre homotopy }
 
 Let us recall the notion of a \textit{clearing map}, also known as \textit{Segre homotopy} (see \cite{E} chapter 1.1 and also \cite{Ca}).
 Let $R$ be a a commutative ring.
 We start with a map $F \in R[X]^n$ of the form (\ref{xh}). Then we can see $F$ as a following sum $F=F_{(1)}+F_{(2)}+ \ldots$, where $F_{(i)}$ is homogeneous of degree $i$.
Following the idea of Segre (see \cite{Ca}) one may instead of $F$ consider a map 
\[\widehat{F}(X) = t^{-1}F(tX)=X+t^{-1}H(tX)=X+\widehat {H}(X).\] 
Here $t$ is a new variable and $\widehat{F} \in (R[t])[X]^n$.
Of course for a two given maps $F$ and $L$ we have 
\[ \widehat{F}(X) \circ \widehat{L}(X) = t^{-1} F(tt^{-1}L(tX))=t^{-1}F(L(tX))=\widehat{F \circ L}(X).\]
Moreover if $G$ is an inverse of $F$, then $\widehat{G}$ is the inverse of $\widehat{F}$. Indeed
\begin{equation}
(\widehat{F} \circ \widehat{G})(X)= \widehat{F}(t^{-1}G(tX))=t^{-1}F(tt^{-1}G(tX))=t^{-1}F(G(tX))=\widehat{F \circ G}(X)=X.
 \label{inv}
\end{equation}
One can check that \[\det(J_{\widehat{F}})(X)=\det{(J_F)}(tX).\]

As mentioned before the map $\widehat{F}$ associated with $F$ is often referred as a clearing map. Let us choose $r \in R$ and  define a new map $\leftidx{^r}{F}$ given by
\[\leftidx{^r}{F}:=\widehat{F}|_{t=r}.\]
So $\leftidx{^r}{F}\in R[X]^n$. The following observation (see \cite{E}, Proposition 1.1.23) justifies the name clearing map.

\begin{lem}\label{dclem}
 Let $R$ be a domain and $K=Fr(R)$ its field of fractions. Let $F \in K[X]^n$ such that $F(0)=0$, $F_{(1)} \in R[X]^n$ and $\det{J_F} \in R^*$, where $R^*$ is the group of units of $R$. Then there exists nonzero $r \in R$ such that $\leftidx{^r}{F}\in R[X]^n$ and 
 $\det{ J_{ ^rF}} \in R^*$.
\end{lem}

To prove this it is enough to choose $r \in R, r \neq 0$ such that for all $i>1$ we have $r \cdot F_{(i)} \in R[X]^n$.
Moreover $\det{ J_{ ^rF}}(X) = \det{J_{\widehat{F}|_{t=r}}}(X)=\det{J_F}(rX)\in R^*$.

\subsection{Algorithm and Segre homotopy}

In this section we discuss behaviour of our algorithm while using Segre homotopy.
We can apply algorithm to both $F$ and $\widehat{F}$. We get two families of polynomial mappings. We establish the notation in the list below.

\[ \begin{array}{c|c}
F&\widehat{F}\\
\hline\\
P_0(X)=X & Q_0(X)=X\\
P_1(X)=P_0(F)-P_0(X)=H(X) & Q_1(X)=Q_0(\widehat{F})-Q_0(X)=\widehat{H}(X)\\
P_2(X)=P_1(F)-P_1(X)=H(F)-H(X) & Q_2(X)=Q_1(\widehat{F})-Q_1(X)=\widehat{H}(\widehat{F})-\widehat{H}(X)\\
  \ldots& \ldots\\
  P_{k+1}(X)=P_k(F)-P_k(X)=(P_k \circ F-P_k)(X)& Q_{k+1}(X)=Q_k(\widehat{F})-Q_k(X)=(Q_k \circ \widehat{F}-Q_k)(X)\\
  \\
  \sum_{j=0}^{m-1}(-1)^j P^i_j(X)=G_i(X)+R^i_m(X) & \sum_{j=0}^{m-1}(-1)^j Q^i_j(X)=\widehat{G}_i(X)+S^i_m(X)
   \end{array}
\]

\begin{lem}\label{shlem1}
 For every $k \in \mathbb{N}$ we have $Q_k(X)= \widehat{P}_k(X)$, where $\widehat{P}_k(X)=t^{-1}P_k(tX)$.
\end{lem}
 \textit{Proof.} For $t=0,1$ thesis holds. 
 Assume that the thesis holds for a given $k \in \mathbb{N}$, then
 \[ Q_{k+1}(X)=Q_k(\widehat{F})-Q_k(X)=(Q_k \circ \widehat{F}-Q_k)(X)= (\widehat{P_k \circ F}-\widehat{P}_k)(X)=\]
 \[= t^{-1}(P_k \circ F)(tX)-t^{-1}P_k(tX)=t^{-1}[P_k\circ F-P_k](tX)=t^{-1}P_{k+1}(tX).\]
 \begin{flushright}
$\square$
\end{flushright}

 \begin{cor}\label{shcor1}
  $F$ is Pascal finite if and only if $\widehat{F}$ is Pascal finite.
 \end{cor}

Now we claim the following.
 
 \begin{lem}\label{shlem2}
  For every $i =1, \ldots, n$ and $m > \frac{D^{n-1}-d_i}{d-1}+1$ we have 
  \[ S^i_m(X)=\widehat{R}^i_m(X) \quad \mathrm{and} \quad S^i_m(\widehat{F})=(-1)^{m+1}\widehat{P}^i_m(X),\]
  where $\widehat{R}^i_m(X)=t^{-1}R^i_m(tX)$.
 \end{lem}
 \textit{Proof.} In theorem given above 
 $G_i(Y_1, \ldots, Y_n)=\sum_{l=0}^{m-1}(-1)^l\widetilde{P}^i_l(Y_1, \ldots, Y_n)$,
  where $\widetilde{P}^i_l$ is the sum of homogeneous summands of $P^i_l$ of degree $ \leq D^{n-1}$. 
  According to (\ref{inv}) it is clear that the inverse polynomial mapping constructed from theorem given above for the map $\widehat{F}$ is exactly $\widehat{G}$.
  If we denote by $\widetilde{Q}^i_l$ the sum of homogeneous summands of $Q^i_l \in (\mathbb{C}[t])[X]$ of degree $ \leq D^{n-1}$, then we obtain
  \[ \widetilde{Q}^i_l(X) = \widehat{\widetilde{P}^i_l}(X):=t^{-1}\widetilde{P}^i_l(tX),\]
  which is clear since the degree with respect to $X$ is the same for $P_l^i$ and $\widehat{P}^i_l$.
  Moreover $\widehat{G}_i(Y)=\sum_{l=0}^{m-1}(-1)^l\widehat{\widetilde{P}^i_l}(Y)$.
  Then by lemma \ref{shlem1}
 \[ S^i_m(X)= \sum_{j=0}^{m-1}(-1)^j Q^i_j(X)-\widehat{G}_i(X)  =
 \sum_{j=0}^{m-1}(-1)^j t^{-1}P^i_j(tX)-t^{-1}{G}_i(tX) = \]
 \[=t^{-1} \bigg( \sum_{j=0}^{m-1}(-1)^j P^i_j-G_i\bigg)(tX)=t^{-1}R^i_m(tX).\]
 And
 \[ S^i_m(\widehat{F})= \widehat{R}^i_m(\widehat{F})=t^{-1}R^i_m(t\widehat{F})=t^{-1}R^i_m(F(tX))=(-1)^{m+1}t^{-1}P^i_m(tX)=(-1)^{m+1}\widehat{P}^i_m(X).\]
  
  \begin{flushright}
$\square$
\end{flushright}
 
 We conclude that the equivalent condition to invertibility of a polynomial map holds when using Segre homotopy.
Due to corollary \ref{shcor1} we already know that for any $r \in R$ mapping $^rF$ is Pascal Finite if and only if $F$ is Pascal Finite. Of course if $G$ is an inverse of $F$, then $^rG$ is an inverse of $^r F$, i.e. $^r(F^{-1})=(^rF)^{-1}$.

 \section{Reductions modulo prime number}
 
 From now on let $R=\mathbb{Z}$ and $K=\Q$.   We can use denominators clearing procedure described above, so we assume that $F \in \Z[X]^n$.
 By $\mathbb{P}$ we denote the set of prime numbers and by $\overline{F}^p \in \mathbb{F}_p[X]^n$ a map obtained from $F$ by reducing coefficients of each 
 $F_i$ modulo given $p \in \mathbb{P}$. 
 If $F \in \mathbb{Z}[X]^n$ is invertible over $\mathbb{Z}$ then $\overline{F}^p$ is invertible over $\mathbb{F}_p$.
We can apply algorithm to both $F$ and $\overline{F}^p$. We get two families of polynomial mappings. The notation is established below.

\[ \begin{array}{c|c}
F&\overline{F}^p\\
\hline\\
P_0(X)=X & V_0(X)=X\\
P_1(X)=P_0(F)-P_0(X)=H(X) & V_1(X)=V_0(\overline{F}^p)-V_0(X)=\overline{H}^p(X)\\
P_2(X)=P_1(F)-P_1(X)=H(F)-H(X) & V_2(X)=V_1(\overline{F}^p)-V_1(X)=\overline{H}^p(\overline{F}^p)-\overline{H}^p(X)\\
  \ldots& \ldots\\
  P_{k+1}(X)=P_k(F)-P_k(X)=(P_k \circ F-P_k)(X)& V_{k+1}(X)=V_k(\overline{F}^p)-V_k(X)=(V_k \circ \overline{F}^p-V_k)(X)\\
   \end{array}
\]

\begin{lem}\label{redlem1}
 For every $k \in \mathbb{N}$ we have $V_k(X)= \overline{P}^p_k(X)$, where $\overline{P}^p_k$ is a reduction modulo $p$ of $P_k$.
\end{lem}
\textit{Proof.} For $t=0,1$ thesis holds. 
 Assume that the thesis holds for a given $k \in \mathbb{N}$. For any $a,b \in R$ we have $\overline{a \cdot p}^p = \overline{a}^p \cdot \overline{b}^p$, so 
 \[ V_{k+1}(X)=V_k \circ \overline{F}^p(X)-V_k(X)=\overline{P_k}^p \circ \overline{F}^p(x) - \overline{P_k}^p =\overline{P_k \circ F - P_k}^p .\]
 \begin{flushright}
$\square$
\end{flushright}

\begin{cor} 
 If $F$ is Pascal finite, then  $\overline{F}^p$ is Pascal finite (for every prime number $p$)
 and an inverse of $\overline{F}^p$ is exactly reduction modulo $p$ of the inverse of $F$, i.e.
 $\overline{F^{-1}}^p=\big(\overline{F}^p \big)^{-1}$.

\end{cor}

\textit{Proof.} First part follows immediately from lemma \ref{redlem1}, $P_k=0$ implies $V_k=0$.
     Moreover by (\ref{inveq}) and lemma \ref{redlem1}
     we obtain
     \[\overline{F^{-1}}^p = \overline{\sum_{l=0}^{m-1} (-1)^l P_l(X)}^p = \sum_{l=0}^{m-1} (-1)^l \overline{P_l}^p(X)=\sum_{l=0}^{m-1} (-1)^l V_l(X)=\big(\overline{F}^p \big)^{-1},\]
     where $m$ is minimal such that $P_m=0$.
 \begin{flushright}
$\square$
\end{flushright}

 Theorem \ref{thm11} stays valid for every field, so in particular for $\mathbb{F}_p$. If $\overline{F}^p$ is invertible, then for every $i =1, \ldots, n$ and $m > \frac{D^{n-1}-d_i}{d-1}+1$ we have
\[ \sum_{j=0}^{m-1}(-1)^j V^i_j(X)=U_i(X)+W^i_m(X)\]
where $U_i(X)$ is a polynomial of degree $ \leq D^{n-1}$, independent of $m$, and $W^i_m(X)$ is a polynomial satisfying $W^i_m(\overline{F}^p)=(-1)^{m+1}V^i_m(X)$, with lower degree $\geq (m-1)(d-1)+d_i>D^{n-1}$.
Moreover the inverse $U$ of $\overline{F}^p$ is given by
\begin{equation}\label{thminveq}
    U_i(Y_1, \ldots, Y_n)=\sum_{l=0}^{m-1}(-1)^l\widetilde{V}^i_l(Y_1, \ldots, Y_n), \, i=1, \ldots, n,
\end{equation}
  where $\widetilde{V}^i_l$ is the sum of homogeneous summands of $V^i_l$ of degree $ \leq D^{n-1}$ and $m$ is an integer $> \frac{D^{n-1}-d_i}{d-1}+1$.

 Observe that due to (\ref{thminveq}) if $F \in \Q[X]^n$ has coefficients in $\Z$, then the inverse $G$ also has coefficients in $\Z$. Using the notation established in theorem \ref{thm11} we claim the following.

 \begin{lem}\label{redlem2} Let $F \in \Z[X]^n$ be invertible.  For every $i =1, \ldots, n$ the following holds.
 \begin{enumerate}
     \item[a)] $U_i(X)=\overline{G_i}^p$
     \item[b)] For every $m > \frac{D^{n-1}-d_i}{d-1}+1$ we have $W^i_m(X)=\overline{R^i_m}^p(X) $.
 \end{enumerate}
  
 \end{lem}
 \textit{Proof.} 
 The degree with respect to $X$ is the same for $P_l^i$ and $\overline{P^i_l}^p$, so it is clear that $\widetilde{\overline{P_l^i}^p}(Y_1, \ldots, Y_n) = \overline{\widetilde{P_l^i}}^p(Y_1, \ldots, Y_n)$.
 Hence \[U_i(Y_1, \ldots, Y_n)=\sum_{l=0}^{m-1}(-1)^l\widetilde{V}^i_l(Y_1, \ldots, Y_n) \quad  \myeqq \quad \sum_{l=0}^{m-1}(-1)^l\widetilde{\overline{P_l^i}^p}(Y_1, \ldots, Y_n) = \]
 \[ =\sum_{l=0}^{m-1} (-1)^l \overline{ \widetilde{P_l^i} }^p (Y_1, \ldots, Y_n) = \overline{ \sum_{l=0}^{m-1} (-1)^l \widetilde{P_l^i} }^p (Y_1,\ldots, Y_n)= G_i(Y_1, \ldots, Y_n). \]
 
  Then
 \[ W^i_m(X)= \sum_{l=0}^{m-1}(-1)^j V^i_l(X)-\overline{G_i}^p(X)\quad  \myeqq \quad 
 \sum_{l=0}^{m-1}(-1)^l \overline{P^i_l}^p(X)-\overline{G_i}^p(X) =\]
 \[=\overline{\Big(  \sum_{l=0}^{m-1}(-1)^l P^i_l- G_i \Big)}^p(X)=\overline{R^i_m}^p(X)  .\]
 And of course
 \[ \overline{R^i_m(F)}^p =  \overline{R^i_m \circ F}^p =\overline{R^i_m}^p \circ \overline{F}^p = W^i_m(\overline{F}^p).\]
  
  \begin{flushright}
$\square$
\end{flushright}
 
 \subsection{Examples of reductions}
 
Observe that $V_k=0$ does not implies $P_k=0$. So if $F$ is Pascal finite, then the number of steps needed to find an inverse of $\overline{F}^p$ is less or equal to the number of steps needed to find an inverse of $F$.

\begin{ex}
Let us consider the following map over $\Q$.
\[
F : \left\{ \begin{array}{ll}
     F_1  = & X_{1}\\
     F_2  = & -\frac{1}{3} X_{1}^{3} + X_{2} \\
     F_3  = & - X_{1}^{2} X_{2} -  X_{1} X_{2}^{2} -  X_{2}^{3} + X_{3} \\
     F_4  = & - X_{1} X_{2}^{2} -  X_{2}^{3} -  X_{1}^{2} X_{3} -  X_{1} X_{2} X_{3} -  X_{2}^{2} X_{3} -  X_{1} X_{3}^{2} -  X_{2} X_{3}^{2} -  X_{3}^{3} + X_{4} 
\end{array} \right.
\]
We can obtain mapping $^3F$ over $\mathbb{Z}$.
\[
^3F : \left\{ \begin{array}{ll}
     ^3F_1  = & X_{1} \\
     ^3F_2  = & -3 X_{1}^{3} + X_{2} \\
     ^3F_3  = & -9 X_{1}^{2} X_{2} - 9 X_{1} X_{2}^{2} - 9 X_{2}^{3} + X_{3} \\
     ^3F_4  = & -9 X_{1} X_{2}^{2} - 9 X_{2}^{3} - 9 X_{1}^{2} X_{3} - 9 X_{1} X_{2} X_{3} - 9 X_{2}^{2} X_{3} - 9 X_{1} X_{3}^{2} - 9 X_{2} X_{3}^{2} - 9 X_{3}^{3} + X_{4} 
\end{array} \right.
\]
Executing algorithm for mapping $^3F$ we obtain $(P_i)_{i \geq 0}$, where $P_i=(P_i^1, P_i^2,P_i^3, P_i^4)$. We can always perform it componentwise. The fourth coordinates $P^4_i$ are presented in table \ref{tab:p_sequnce}. The algorithm executed for mapping $\overline{^3F}^2$ produces sequence $(V_i)_{i \geq 0}$, where $V_i=(V_i^1, V_i^2,V_i^3, V_i^4)$. The fourth coordinates $V_i^4$ are presented in table \ref{tab:v_sequnce}.
 One can observe that after reduction modulo $p$, number of steps which are necessary to obtain the inverse can decrease.
\end{ex}

\begin{table}
\parbox{.45\linewidth}{
\centering
 \begin{tabular}{|c|c|c|c|}
    \hline
    {\footnotesize Element} & {\footnotesize Number} & {\footnotesize Degree} & {\footnotesize Ldegree} \\
    &{\footnotesize of monomials}&& \\\hline
    $P_{0}^4$ & 1 & 1 & 1 \\
    $P_{1}^4$ & 8 & 3 & 3 \\
    $P_{2}^4$ & 39 & 9 & 5\\
    $P_{3}^4$ & 97 & 27 & 7 \\
    $P_{4}^4$ & 79 & 27 & 9\\
    $P_{5}^4$ & 61 & 27 & 11\\
    $P_{6}^4$ & 46 & 27 & 13\\
    $P_{7}^4$ & 34 & 27 & 15\\
    $P_{8}^4$ & 24 & 27 & 17\\
    $P_{9}^4$ & 16 & 27 & 19\\
    $P_{10}^4$ & 10 & 27 & 21 \\
    $P_{11}^4$ & 6 & 27 & 23\\
    $P_{12}^4$ & 3 & 27 & 25\\
    $P_{13}^4$ & 1 & 27 & 27\\
    $P_{14}^4$ & 0 & & \\ \hline
    \end{tabular}
\caption{Sequence $(P^4_i)$}
    \label{tab:p_sequnce}
}
\hspace{10pt} 
\parbox{.45\linewidth}{
\centering
\begin{tabular}{|c|c|c|c|}
    \hline
    {\footnotesize Element} & {\footnotesize Number} & {\footnotesize Degree} & {\footnotesize Ldegree} \\
    &{\footnotesize of monomials}&&\\ 
    \hline
    $V_{0}^4$ & 1 & 1 & 1 \\
    $V_{1}^4$ & 8 & 3 & 3 \\
    $V_{2}^4$ & 22 & 9 & 5\\
    $V_{3}^4$ & 36 & 27 & 7 \\
    $V_{4}^4$ & 17 & 21 & 9\\
    $V_{5}^4$ & 0 &  & \\ \hline
    \end{tabular}
    \caption{Sequence $(V^4_i)$}
    \label{tab:v_sequnce}
}
\end{table}

One can ask about the property of \emph{not} being Pascal finite. What happens when we reduce the coefficients modulo prime number? Does the property holds? An example given below answers this question.

\begin{ex}
Let us consider the following map over $\Q$, which is a representative of the eighth class in Hubbers classification of cubic homogeneous polynomial maps over fields of characteristic zero in dimension 4 (see \cite{E}, Theorem 7.1.2). $F$ is not Pascal finite (see \cite{ABCH2}, Remark 3.2).
\[
F : \left\{ \begin{array}{ll}
     F_1  = & X_{1}\\
     F_2  = & -1/3 X_1^3 + X_2 \\
     F_3  = & -X_1^2 X_2 - 7 X_1 X_2^2 - 7 X_2^3 + 7 X_1 X_2 X_3 + 7 X_2^2 X_3 + 49 X_2^2 X_4 + X_3 \\
     F_4  = &  -7 X_1 X_2^2 - 7 X_2^3 - X_1^2 X_3 - 2 X_1 X_2 X_3 - X_2^2 X_3 - 7 X_1 X_2 X_4 - 7 X_2^2 X_4 + X_4
\end{array} \right.
\]

We consider $^3F \in \mathbb{Z}[X]^4$, which by corollary \ref{shcor1} is not Pascal finite.
\[
^3F : \left\{ \begin{array}{ll}
     ^3 F_1  = & X_{1} \\
     ^3 F_2  = & -3 X_1^3 + X_2 \\
     ^3 F_3  = & -9 X_1^2 X_2 - 63 X_1 X_2^2 - 63 X_2^3 + 63 X_1 X_2 X_3 + 63 X_2^2 X_3 + 441 X_2^2 X_4 + X_3 \\
     ^3 F_4  = & -63 X_1 X_2^2 - 63 X_2^3 - 9 X_1^2 X_3 - 18 X_1 X_2 X_3 - 9 X_2^2 X_3 - 63 X_1 X_2 X_4 - 63 X_2^2 X_4 + X_4
\end{array} \right.
\]

Now we reduce all coefficients of ${^3F}$ modulo $5$.
\[
\overline{^3F}^5 : \left\{ \begin{array}{l}
     \overline{^3F}^5_1 =   X_1 \\
     \overline{^3F}^5_2 = 2 X_1^3 + X_2\\
     \overline{^3F}^5_3 = X_1^2 X_2 + 2 X_1 X_2^2 + 2 X_2^3 - 2 X_1 X_2 X_3 - 2 X_2^2 X_3 + X_2^2 X_4 + X_3 \\
     \overline{^3F}^5_4 = 2 X_1 X_2^2 + 2 X_2^3 + X_1^2 X_3 + 2 X_1 X_2 X_3 + X_2^2 X_3 + 2 X_1 X_2 X_4 + 2 X_2^2 X_4 + X4
\end{array} \right. 
\]

The algorithm executed for mapping $\overline{^3F}^5$ produces sequence $(V_i)$. The fourth coordinates of its elements are presented in table \ref{tab:hub8par7red5_v_sequnce}. Observe that $\overline{^3F}^2$ is not Pascal finite. It can be proved that the lower degree of $V^4_k$ is exactly $2k+1$. 

\begin{table}[h]
    \centering
    \begin{tabular}{|c|c|c|c|}
    \hline
    {\footnotesize Element} & {\footnotesize Number} & {\footnotesize Degree} & {\footnotesize Ldegree} \\
    &{\footnotesize of monomials}&& \\\hline
    $V_{0}^4$ & 1 & 1 & 1 \\
    $V_{1}^4$ & 7 & 3 & 3 \\
    $V_{2}^4$ & 27 & 9 & 5\\
    $V_{3}^4$ & 40 & 15 & 7 \\
    $V_{4}^4$ & 50 & 19 & 9\\
    $V_{5}^4$ & 61 & 23 & 11\\
    $V_{6}^4$ & 71 & 27 & 13\\
    $V_{7}^4$ & 82 & 31 & 15\\
    $V_{8}^4$ & 92 & 35 & 17\\
    $V_{9}^4$ & 103 & 39 & 19\\
    $V_{10}^4$ & 113 & 43 & 21\\
    $V_{11}^4$ & 124 & 47 & 23\\
    $V_{12}^4$ & 134 & 51 & 25\\
    $V_{13}^4$ & 145 & 55 & 27\\
    $V_{14}^4$ & 155 & 59 & 29 \\ \hline
    \end{tabular}
    \caption{Sequence $(V^4_i)$}
    \label{tab:hub8par7red5_v_sequnce}
\end{table}

Let us now reduce $^3F$ modulo $7$.
\[
\overline{^3F}^7 : \left\{ \begin{array}{l}
     \overline{^3F}^7_1 =   X_1 \\
     \overline{^3F}^7_2 = -3 X_1^3 + X_2\\
     \overline{^3F}^7_3 = -2 X_1^2 X_2 + X_3 \\
     \overline{^3F}^7_4 = -2 X_1^2 X_3 + 3 X_1 X_2 X_3 - 2 X_2^2 X_3 + X_4
\end{array} \right. 
\]
$\overline{^3F}^7$ is Pascal finite since it is triangular (see \cite{ABCH2}, Corollary 2.1.).
We conclude that reduction modulo prime number of a given not Pascal finite map can be both Pascal finite or not Pascal finite, depending on the choice of a prime number $p$.


\end{ex}

 \section{Finding an inverse of polynomial map with integer coefficients}

 Here arises a question, can we somehow retrieve $F^{-1}$ for $F\in \mathbb{Z}[X]^n$ knowing $\big(\overline{F}^p \big)^{-1}$ for $p \in S$, where $S$ is some finite
 subset of the set $\mathbb{P}$ of all prime numbers?

 \subsection{An introductory example}
 Consider  $F: \Q^4 \rightarrow \Q^4$ given by
\[  F: \left\{ \begin{array}{l}
 F_1=X_1\\
 F_2=X_2+X_3X_4^2+X_1X_2X_4-X_1X_4^2+X_2X_4^2 -X_4^3+\frac{1}{3}X_1^3\\
 F_3=X_3-X_1X_3X_4-X_1X_4^2-X_3X_4^2-X_4^3-X_1^2X_2-2X_1X_2X_4-X_2X_4^2\\
 F_4=X_4
\end{array} \right. .\]
We clear denominators and obtain the following map $^3 F \in \mathbb{Z}[X]^4$.
\begin{equation}
  ^3 F :\left\{ \begin{array}{l}
    {^3 F_1} = X_1  \\
    {^3F_2}  =X_2 -3 X_{1}^{3} + 9 X_{1} X_{2} X_{4} - 9 X_{1} X_{4}^{2} + 9 X_{2} X_{4}^{2} + 9 X_{3} X_{4}^{2} - 9 X_{4}^{3} \\
    {^3F_3} = X_3 -9 X_{1}^{2} X_{2} - 18 X_{1} X_{2} X_{4} - 9 X_{1} X_{3} X_{4} - 9 X_{1} X_{4}^{2} - 9 X_{2} X_{4}^{2} - 9 X_{3} X_{4}^{2} - 9 X_{4}^{3}  \\
    {^3F_4} = X_4 
\end{array}  \right.  
\label{introex}
\end{equation}
We can find its inverse $^3G$ using the algorithm.
\[
^3 G : \left\{ \begin{array}{ll}
     ^3 G_1   = & X_{1} \\
     ^3 G_2   = &X_2 -27 X_{1}^{4} X_{4} - 27 X_{1}^{3} X_{4}^{2} - 81 X_{1}^{2} X_{4}^{3} - 243 X_{1} X_{4}^{4} - 162 X_{4}^{5} + 3 X_{1}^{3} - 9 X_{1} X_{2} X_{4} \\ & + 9 X_{1} X_{4}^{2} - 9 X_{2} X_{4}^{2} - 9 X_{3} X_{4}^{2} + 9 X_{4}^{3} \\
    ^3 G_3  = &X_3 + 27 X_{1}^{5} + 54 X_{1}^{4} X_{4} + 108 X_{1}^{3} X_{4}^{2} + 324 X_{1}^{2} X_{4}^{3} + 405 X_{1} X_{4}^{4} + 162 X_{4}^{5} \\ & + 9 X_{1}^{2} X_{2} + 18 X_{1} X_{2} X_{4} + 9 X_{1} X_{3} X_{4} + 9 X_{1} X_{4}^{2} + 9 X_{2} X_{4}^{2} + 9 X_{3} X_{4}^{2} + 9 X_{4}^{3} \\
    ^3 G_4  = &X_{4} 
\end{array} \right.
\]
We reduce coefficients of ${^3F}$ modulo $5$ and obtain $\overline{^3F}^5  \in \mathbb{F}_5[X]^4$.
\begin{equation}
\overline{^3F}^5 : \left\{ \begin{array}{l}
     \overline{^3F}^5_1 =   X_1 \\
     \overline{^3F}^5_2 = X_2+2 X_{1}^{3} -  X_{1} X_{2} X_{4} + X_{1} X_{4}^{2} -  X_{2} X_{4}^{2} -  X_{3} X_{4}^{2} + X_{4}^{3} \\
     \overline{^3F}^5_3 = X_3+X_{1}^{2} X_{2} + 2 X_{1} X_{2} X_{4} + X_{1} X_{3} X_{4} + X_{1} X_{4}^{2} + X_{2} X_{4}^{2} + X_{3} X_{4}^{2} + X_{4}^{3} \\
     \overline{^3F}^5_4 = X_4
\end{array} \right. 
\label{fmod5}
\end{equation}
Using our algorithm we can find its inverse $\overline{^3G}^5$ over $\F_5$.
\begin{equation}
\overline{^3G}^5 : \left\{ \begin{array}{ll}
     \overline{^3G}^5_1 = &   X_1 \\
     \overline{^3G}^5_2 = & X_2+3 X_{1}^{4} X_{4} + 3 X_{1}^{3} X_{4}^{2} -  X_{1}^{2} X_{4}^{3} + 2 X_{1} X_{4}^{4} + 3 X_{4}^{5} + 3 X_{1}^{3}  \\ & + X_{1} X_{2} X_{4} -  X_{1} X_{4}^{2} + X_{2} X_{4}^{2} + X_{3} X_{4}^{2} -  X_{4}^{3}  \\
     \overline{^3G}^5_3 = & X_3+2 X_{1}^{5} -  X_{1}^{4} X_{4} + 3 X_{1}^{3} X_{4}^{2} -  X_{1}^{2} X_{4}^{3} + 2 X_{4}^{5} -  X_{1}^{2} X_{2} \\ & + 3 X_{1} X_{2} X_{4} -  X_{1} X_{3} X_{4} -  X_{1} X_{4}^{2} -  X_{2} X_{4}^{2} -  X_{3} X_{4}^{2} -  X_{4}^{3}  \\
     \overline{^3G}^5_4 = & X_4
\end{array} \right. 
\label{gmod5}
\end{equation}
As one can see  $\overline{^3F}^5$ is invertible over $\F_5$, hence ${^3F}$ is invertible over $\Z$ and $F$ is invertible over $\mathbb{Q}$.
Observe at this point that in formulas (\ref{fmod5}) and (\ref{gmod5}) we have some freedom of choosing a representative of a given congruence class. However we decide to always choose the one with the smallest absolute value. For example we see $\mathbb{F}_5=\{-2,-1,0,1,2\}$ instead of $\{0,1,2,3,4\}$. In this way we can deal with negative coefficients. We comment on this choice in the next section.

Now one can ask if it is possible to retrieve the inverse $^3G$ of $^3F$ knowing $\overline{^3G}^5$. This information is clearly not enough, however we can find such inverses $\overline{^3G}^p = (\overline{^3G_1}^p, \overline{^3G_2}^p,\overline{^3G_3}^p, \overline{^3G_4}^p)$ over $\F_p$, for $p \in S$, where $S \subset \mathbb{P}$ is finite.  We can consider it componentwise. We distinguish monomials appearing in  $\overline{^3G_i}^p$ and consider sequences of coefficients appearing alongside each monomial.  We present coefficients appearing in the second coordinate of the inverse mappings in table \ref{t:1}.
\begin{table}[h]
    \centering
    \begin{tabular}{|c|r||r|r|r|r|r|r|r|r|r||c|c|c|}
    \hline
        \footnotesize{Monomial} & $\mathbf{^3G}$&$\overline{^3G}^5$ & $\overline{^3G}^7$ & $\overline{^3G}^{11}$ & $\overline{^3G}^{13}$ &
        $\overline{^3G}^{17}$ &
        $\overline{^3G}^{19}$ &
        $\overline{^3G}^{23}$ &
        $\overline{^3G}^{59}$ &
        $\overline{^3G}^{61}$ &
        $\scriptstyle{N_1=385}$ &
        $\scriptstyle{N_2=5005}$ & $\scriptstyle{N_3=85085 }$ \\ \hline 
         $X_4^3$ &9& -1 & 2 & -2 & -4 &-8&\textbf{9}&\textbf{9}&\textbf{9}&\textbf{9}& \textbf{9} & \textbf{9} &\textbf{9}\\ \hline 
         $X_4^5$ &-162& -2 & -1 & 3 & -6 &8&9&-1&15&21& \textbf{-162} &\textbf{ -162} &\textbf{-162}\\ \hline
         $X_3X_4^2$ &-9& 1 & -2 & 2 & 4 &8&\textbf{-9}&\textbf{-9}&\textbf{-9}&\textbf{-9}& \textbf{-9} & \textbf{-9} &\textbf{-9}\\ \hline
         $X_2$ &1& 1 & 1 & 1 & 1 &1&1&1&1&1& 1 & 1 &1\\ \hline
         $X_2X_4^2$ &-9& 1 & -2 & 2 & 4 &8&\textbf{-9}&\textbf{-9}&\textbf{-9}&\textbf{-9}& \textbf{-9} &\textbf{ -9} &\textbf{-9} \\ \hline
         $X_1X_4^2$ &9& -1 & 2 & -2 & -4 &-8&\textbf{9}&\textbf{9}&\textbf{9}&\textbf{9}& \textbf{9} & \textbf{9} &\textbf{9}\\ \hline
         $X_1X_4^4$ &-243& 2 & 2 & -1 & 4 &-5&4&10&-7&1& 142 & \textbf{-243}&\textbf{-243} \\ \hline
         $X_1X_2X_4$ &-9& 1 & -2 & 2 & 4 &8&\textbf{-9}&\textbf{-9}&\textbf{-9}&\textbf{-9}& \textbf{-9} & \textbf{-9} &\textbf{-9} \\ \hline
         $X_1^2X_4^3$ &-81& -1 & 3 & -4 & -3 &4&-5&11&-22&-20& \textbf{-81} &\textbf{ -81}&\textbf{-81} \\ \hline
         $X_1^3$ &3& -2 & \textbf{3} & \textbf{3}  & \textbf{3}  &\textbf{3} &\textbf{3} &\textbf{3} &\textbf{3} &\textbf{3} & \textbf{3} & \textbf{3} &\textbf{3}\\ \hline
         $X_1^3X_4^2$ &-27& -2 & 1 & -5 & -1 &7&-8&-4&\textbf{-27}&\textbf{-27}& \textbf{-27} &\textbf{ -27} &\textbf{-27}\\ \hline
         $X_1^4X_4$ &-27& -2 & 1 & -5 & -1 &7&-8&-4&\textbf{-27}&\textbf{-27}& \textbf{-27} & \textbf{-27} &\textbf{-27}\\ \hline
    \end{tabular}
    \caption{Coefficients of  $\overline{^3G}^{p}_2$ for various $p$  (representatives with the smallest absolute value).}
    \label{t:1}
\end{table}

We observe stabilization of coefficients in all but three rows of the table \ref{t:1}.  One can suspect that after considering $p$ large enough one can be able to obtain stabilization also in the three remaining rows. 
Instead of investigating
many prime numbers we use Chinese Remainder Theorem (see for example \cite{IR}, chapter 3) which allows us to get an element  of a ring $ \mathbb{Z}/_{N\mathbb{Z}}$ for relatively big $N$.
Denote $N_1 = 5 \cdot 7 \cdot 11 = 385$, $N_2 = 5 \cdot 7 \cdot 11 \cdot 13 =5005 $ and $N_3 = 5 \cdot 7 \cdot 11 \cdot 13 \cdot 17 =85085 $. Values in the column $N_1$ are coefficients in the ring $\mathbb{Z}/_{N_1 \mathbb{Z}}$, calculated by the Chinese Remainder Theorem for moduli $5,7,11$. Similarly for $N_2$ and $N_3$. Now indeed we can observe stabilization of coefficients in all rows. For example coefficient of $X_4^3$ appearing in $^3G_2 \in \mathbb{Z}[X]$ is congruent to 9 modulo 385 and modulo 5005 and  modulo 85085. Let us assume then that this coefficient is equal to 9. We repeat such a procedure for every monomial and we obtain the following polynomial
\[ \begin{array}{ll}
    T_2 = & X_2-27 X_{1}^{4} X_{4} - 27 X_{1}^{3} X_{4}^{2} - 81 X_{1}^{2} X_{4}^{3} - 243 X_{1} X_{4}^{4} - 162 X_{4}^{5} + 3 X_{1}^{3} - 9 X_{1} X_{2} X_{4} \\
     & + 9 X_{1} X_{4}^{2} - 9 X_{2} X_{4}^{2} - 9 X_{3} X_{4}^{2} + 9 X_{4}^{3} 
\end{array} \]
One can check that $ T_2={}^3G_2$.
This allows us to suspect, that some algorithmic method for choosing particular coefficients while retrieving $^3G$ can be proposed.

\subsection{Stabilization of coefficients while reducing modulo prime number}

Let $F \in \mathbb{Z}[X]^n$ be a polynomial automorphism of the form (\ref{xh}).  Let us choose a finite subset of primes $S \subset \mathbb{P}$. Denote by $(\overline{F}^p)_{p \in S}$ a sequence of reductions of $F$ modulo prime numbers $p$. Our goal is to retrieve its inverse $G$ by considering  sequence of inverse maps $(\overline{G}^p)_{p \in S}$ obtained by performing the algorithm for each $\overline{F}^p$. Here $\overline{G}^p = (\overline{G_1}^p, \ldots, \overline{G_n}^p)$. We can consider it componentwise, each $\overline{G_i}^p$ separately. We distinguish monomials appearing in  $(\overline{G_i}^p)_{p \in S}$ and consider sequences of coefficients appearing alongside each monomial, i.e. a alongside each product of the form $X_1^{a_1}X_2^{a_2}\cdot \ldots \cdot X_n^{a_n}$, where $a_i \in \{0,1,2, \ldots\}$ for $i  \in \{1,\ldots, n\}$. 
If $M$ is a monomial appearing in  $(\overline{G_i}^p)_{p \in S}$, then we obtain a finite sequence of coefficients $(\overline{\alpha_M}^p) _{p \in S}$. Here we understand $\overline{\alpha_M}^p$ as a representative of congruence class in $\mathbb{Z}/p\mathbb{Z}$ with the smallest absolute value, i.e. an element of $\mathbb{Z}$.

\begin{df}
We say that the coefficient of a monomial $M$ \emph{stabilizes} when there exists $p_0 \in \mathbb{P}$ such that for every $p \in \Z$, $p\geq p_0$ we have $\overline{\alpha_M}^p=\overline{\alpha_M}^{p_0}$.
\end{df}

Observe that if $\alpha \in \mathbb{Z}$ is a coefficient of a monomial in $G$, then for every $p \in \mathbb{Z}$ (not necessarily prime) the following holds

\begin{equation}
    p>|2\alpha | \quad \Rightarrow \quad\overline{\alpha }^p=\alpha .
    \label{sav}
\end{equation} So when we are performing reductions modulo consecutive prime numbers, then the coefficient appearing in each row of the table \ref{t:1} will finally stabilize, irrespective of the sign of $\alpha $, since we decided to always choose a representative with the smallest absolute value.

Here arise two questions about proposed way of treating the problem.
By lemma \ref{redlem2} monomials appearing in $\overline{G_i}^p$ are those appearing in $G_i \in \mathbb{Z}[X]$ (maybe some of them with zero coefficient).
A priori we do not know $G$, so we consider monomials appearing in at least one of $\overline{G_i}^p$. 
One can ask, how to check, that when performing reductions modulo some finite set of prime numbers we obtain all monomials of $G$.

\begin{ex}
Consider $\alpha = 255255$. Since $\alpha =3 \cdot 5 \cdot 7 \cdot 11 \cdot 13 \cdot 17$, then
\[ \overline{\alpha }^3=\overline{\alpha }^5=\overline{\alpha }^7=\overline{\alpha }^{11}=\overline{\alpha }^{13}=\overline{\alpha }^{17}=0. \]
However $\overline{\alpha }^{p} \neq 0$, for every prime number $ p >17$.
\label{zerocoeff}
\end{ex}

 Another question is, when we actually observe a stabilization? When one can be sure that if $ \overline{\alpha } ^p = \alpha$, then for every $q \in \mathbb{Z}$, $q>p$, we have $\overline{\alpha }^q = \alpha$?

\begin{ex}
Consider $\alpha = 255257$. Since $\alpha =2+3 \cdot 5 \cdot 7 \cdot 11 \cdot 13 \cdot 17$, then
\[ \overline{\alpha }^3=\overline{\alpha }^5=\overline{\alpha }^7=\overline{\alpha }^{11}=\overline{\alpha }^{13}=\overline{\alpha }^{17}=2. \]
However $\overline{\alpha }^{19}=11$, $\overline{\alpha }^{23}=3$ etc. By (\ref{sav}) we have $\overline{\alpha }^{N}=\alpha $ for every $N>2\alpha =510514$.
\label{stabilization}
\end{ex}

 Examples \ref{zerocoeff} and \ref{stabilization}  illustrate two problems appearing during retrieving coefficients of $G$.
The input of the algorithm is a polynomial automorphism $F$.  But we do not know anything about coefficients of the inverse mapping $G$. If we would be able do determine the coefficient of $G$ with the biggest absolute value, then by (\ref{sav})  we would know when we can actually observe stabilization. 
However investigation of to many prime numbers or performing reduction modulo big prime number will not allow us to decrease the amount of time needed. 
The idea is to use Chinese Remainder Theorem for a given finite subset of primes to find an element $\overline{\alpha }^N$  of a ring $ \mathbb{Z}/_{N\mathbb{Z}}$ for relatively big $N$ in order to confirm, that we actually observe a stabilization. 
Also a decision procedure to answer if obtained set of monomials is the whole set of monomials of $G$ is needed.

 \subsection{Estimation of the coefficients of the inverse map}

For a polynomial $T(X_1, \ldots X_n)$ over an arbitrary field we can determine the number of monomials appearing in $T$. Let us denote it by $l(T)$ and call it the \emph{length} of polynomial $T$. If $T=(T_1, \ldots, T_n)$ is a polynomial mapping, then we set $l(T) = \max \{l(T_1), \ldots, l(T_n)\}$.

For once given $F=(F_1, \ldots, F_n) \in \mathbb{Z}[X]^n$ of the form $(\ref{xh})$ we know its degree $D$, lower degree $d$ of the map $H=(H_1, \ldots, H_n)$, number of variables $n$ and we can determine its length $l(F) = \max _{i=1, \ldots ,n} l(F_i)$. Let $Z_F$ denote the set of all coefficients of monomials appearing in $F$ and $Z_G$ denote the set of all coefficients of monomials appearing in $G$. Let $B=\max \{ |\alpha |: \alpha  \in Z_F\}$ and $A=\max \{ |\alpha |: \alpha  \in Z_G\}$. We would like to find an upper bound for $A$ depending only on $D,d, n, l(F)$ and $B$.

In order to estimate $A$ we perform the algorithm for $F$. We consider each polynomial map $P_k=(P_k^1, \ldots , P_k^n)$ and estimate its length $l(P_k)=\max_{i=1, \ldots, n} l(P_k^i)$. 
By \cite{ABCH} lemma 2.2 we know that $\deg(P_k) \leq D^k$. If $Z_{P_k^i}$ is the set of coefficients of monomials appearing in $P_k^i$, then we set $B_{ki}=\max\{|z|: z \in Z_{P_k^i}\}$ and $B_k=\max_{i=1,\ldots, n} B_{ki}$. We start with the following.

\[ \begin{array}{l|cc}
\mathrm{polynomial \,\,map}&\mathrm{length} & \mathrm{coefficient}\\
\hline\\
P_0 = Id& l(P_0)=1& B_0=1\\
\\
P_1=F-Id& l(P_1) = l(F)-1 & B_1=B\\
   \end{array}
\]

 \begin{lem}
  Let $F=(F_1, \ldots, F_n) \in \mathbb{Z}[X]^n$ be a polynomial map of the form (\ref{xh}) of degree $D$ and let $\{P_k\} _{k \geq 0}$ be a sequence of polynomial mappings obtained when performing an algorithm for $F$. Then for every $k=1,2,\ldots$ we have
  \begin{equation}
      l(P_{k+1}) \leq l(P_k) \cdot [l(F)^{D^k}+1].
      \label{stlem1eq}
  \end{equation}
  \label{stlem1}
 \end{lem}
 
  \textit{Proof.} 
 For $k = 1$ we
  consider $P_2=P_1 \circ F - P_1 = P_1 \circ (F - Id)$. Observe that $P_1 \circ F$ has exactly $l(P_1)$ monomials when seen as a polynomial map in variable $F$ and at most $l(P_1) \cdot[l(F)^{\deg{P_1}}]$ monomials when seen as a polynomial map in variable $X$. So \[l(P_2) \leq l(P_1) \cdot[l(F)^{\deg{P_1}}] + l(P_1) = l(P_1) \cdot [l(F)^{\deg{P_1}}+1]. \]

 Let us assume
 that the thesis holds for some $k>1$. Then $P_{k+1} = P_k \circ F - P_k$. 
 Observe that $P_k \circ F$ has exactly $l(P_k)$ monomials when seen as a polynomial map in variable $F$ and at most $l(P_k) \cdot[l(F)^{\deg{P_k}}]$ monomials when seen as a polynomial map in variable $X$. Since $\deg{P_k} \leq D^k$, we get the thesis.
  \begin{flushright}
$\square$
\end{flushright}

\begin{cor}
Let $F$ be as above. Then for every $k=1,2,\ldots$ we have
  \begin{equation}
      l(P_{k+1}) \leq \big(l(F)-1 \big) \cdot \prod _{j=1}^{k} \big[ l(F)^{D^j} +1 \big] .
      \label{stcor1eq}
  \end{equation} 
  \label{stcor1}
\end{cor}

\textit{Proof.} 
 The thesis holds for $k = 1$. Let us assume
 that the thesis holds for some $k>1$, i.e. $l(P_{k}) \leq \big(l(F)-1 \big) \cdot \prod _{j=1}^{k-1} \big[ l(F)^{D^j} +1 \big] 
 .$ Then by lemma \ref{stlem1} we get 
 \[ l(P_{k+1}) \leq l(P_k) \cdot [l(F)^{D^k}+1]=\big(l(F)-1 \big) \cdot \prod _{j=1}^{k} \big[ l(F)^{D^j} +1 \big] .\]

  \begin{flushright}
$\square$
\end{flushright}

Let us denote an obtained upper bound for $l(P_k)$ by $l_k$, i.e.  
\[l_k = \big(l(F)-1 \big) \cdot \prod _{j=1}^{k-1} \big[ l(F)^{D^j} +1 \big].\]
The sequence $(l(P_k))_{k \geq 0}$ does not have to be increasing, but the sequence $(l_k)_{k\geq0}$ is always increasing.
Let us now estimate elements of the sequence $(B_k)_{k\geq 0}$. We will give an upper bound in worst possible case. So we assume that $\pm B$ appears in monomial of $F$ of degree $D$.

\begin{lem}
Let $F$ be as above. Then for every $k=1,2,\ldots$ we have
\begin{equation}
    B_{k+1} \leq  B_k \cdot B^{D^k} \cdot l_{k+1}.
 \label{stlem2eq}   
\end{equation}
\label{stlem2}
\end{lem}

\textit{Proof.} 
 If $k = 1$, then $B_2 \leq B \cdot B^D \cdot \beta$, where $\beta$ is a number coming from addition or substraction of monomials in $P_1 \circ F-P_1$. 
Hence $\beta \leq l_2$ and $B_2 \leq B^{1+D} \cdot l_2$.
 
 Let us assume
 that the thesis holds for some $k>1$.
 We have $P_{k+1}=P_k \circ F - P_k$. The module of a coefficient with the largest module in $P_k \circ F$ is less or equal to $B_k \cdot B^{\deg{P_k}} \cdot \gamma$, where $\gamma$ is a number coming from addition or substraction of monomials in $P_k \circ F - P_k$.
Hence $\gamma \leq l_{k+1}$ and $B_{k+1} \leq B_k \cdot B^{D^k} \cdot l_{k+1}$.
  \begin{flushright}
$\square$
\end{flushright}

\begin{cor}
Let $F$ be as above. Then for every $k=1,2,\ldots$ we have
\begin{equation}
B_{k+1} \leq B^{1+D+D^2+\ldots + D^k} \cdot \big(l(F)-1 \big)^k \cdot \prod _{j=1}^{k} \Big( l(F)^{D^j} + 1\Big)^{k+1-j}.
\label{stcor2eq}
\end{equation}
\label{stcor2}
\end{cor}

\textit{Proof.} 
Since $B_2 \leq B^{1+D} \cdot l_2 = B^{1+D} \cdot (l(F)-1)(l(F)^D+1)$, then the thesis holds for $k=1$. Let us assume
 that the thesis holds for some $k>1$. By lemma \ref{stlem2} we have 
\[B_{k+1} \leq B^{1+D+D^2+\ldots + D^k} \cdot \prod _{j=2}^{k+1}l_j.\]
Indeed, if $B_k \leq B^{1+D+D^2+\ldots + D^{k-1}} \cdot \prod _{j=2}^{k}l_j$, then $B_{k+1} \leq B_k \cdot B^{D^k} \cdot l_{k+1} = B^{1+D+D^2+\ldots + D^k} \cdot \prod _{j=2}^{k+1}l_j$. Moreover by (\ref{stcor1eq}) we obtain
\[ \prod _{j=2}^{k+1}l_j =  \prod _{j=2}^{k+1} \Big[ \, \big(l(F)-1 \big)  \prod _{s=1}^{j-1} \big[ l(F)^{D^s} +1 \big] \, \Big] = \big(l(F)-1 \big)^k \cdot \prod _{j=1}^{k} \Big( l(F)^{D^j} + 1\Big)^{k+1-j}.\]
Hence we get the thesis.
  \begin{flushright}
$\square$
\end{flushright}

As mentioned before we compute a bound for $B_k$ in the the worst possible case. Let us denote an obtained bound by $b_k$, i.e.  
\[b_k:=B^{1+D+D^2+\ldots + D^{k-1}} \cdot \big(l(F)-1 \big)^{k-1} \cdot \prod _{j=1}^{k-1} \Big( l(F)^{D^j} + 1\Big)^{k-j} .\] Observe that the sequence $(B_k)_{k \geq 0}$ does not have to be increasing, but the sequence $(b_k)_{k\geq0}$ is always increasing.

\begin{thm}
Let $F \in \mathbb{Z}[X]^n$ be a polynomial automorphism of the form (\ref{xh}) with the inverse $G\in \mathbb{Z}[X]^n$. 
Let $Z_F$ denote the set of all coefficients of monomials appearing in $F$ and $Z_G$ denote the set of all coefficients of monomials appearing in $G$. Let $B=\max \{ |z|: z \in Z_F\}$ and denote $A=\max \{ |z|: z \in Z_G\}$. Then 
\begin{equation}
    A \leq B^{\sum_{i=0}^{\mu -2}D^{i}} \cdot (l(F)-1)^{\mu - 2} \prod_{j=1}^{\mu-2}\Big( l(F)^{D^j}+1\Big)^{\mu -1-i} \cdot \Big[ l(F) + (l(F)-1)\sum_{r=2}^{\mu-1} \prod_{s=1}^{r-1}(l(F)^{D^s} +1)\Big],
    \label{estim}
\end{equation}
where $\mu := \lfloor \frac{D^{n-1}-d}{d-1}+1 \rfloor +1$.
\label{stthm}
\end{thm}

\textit{Proof.} 
Let $\{P_k\} _{k \geq 0}$ be a sequence of polynomial mappings obtained when performing an algorithm for $F$.
By theorem \ref{symthm} the inverse $G$ of $F$ is given by  $ G=\sum_{i=0}^{\mu-1}(-1)^l\widetilde{P}_i $. Of course $l(\widetilde{P}_i) \leq l(P_i) \leq l_i$. Hence
$A \leq b_{\mu-1} \cdot  l(G) \leq b_{\mu-1} \cdot \sum_{i=0}^{\mu - 1} l_i $. By (\ref{stcor2eq})
\[b_{\mu-1} = B^{1+D+D^2+\ldots + D^{\mu-2}} \cdot \big(l(F)-1 \big)^{\mu-2} \cdot \prod _{j=1}^{\mu-2} \Big( l(F)^{D^j} + 1\Big)^{\mu-1-j}.\]
Moreover  we have
\[l_i =  \big(l(F)-1 \big) \cdot \prod _{s=1}^{i-1} \big[ l(F)^{D^s} +1 \big] \]
and 
\[\sum_{i=0}^{\mu - 1} l_i =1+(l(F)-1) +  \sum_{i=2}^{\mu - 1} l_i = l(F) +  (l(F)-1)\sum_{r=2}^{\mu-1} \prod_{s=1}^{r-1}(l(F)^{D^s} +1).\]
  \begin{flushright}
$\square$
\end{flushright}

 \subsection{Retrieving the inverse map}
 
 Theorem \ref{stthm} allows us to propose a procedure of retrieving the inverse of a polynomial automorphism with integer coefficients.
We use the notation established in previous sections. 
Let $F \in \mathbb{Z}[X]^n$ be a polynomial automorphism of the form (\ref{xh}) with the inverse $G\in \mathbb{Z}[X]^n$. We choose finite subset $S \subset \mathbb{P}$ and consider sequences $(\overline{F}^p)_{p \in S}$ and $(\overline{G}^p)_{p \in S}$. Here $\overline{G}^p = (\overline{G_1}^p, \ldots, \overline{G_n}^p)$.  We distinguish monomials appearing in  $(\overline{G_i}^p)_{p \in S}$ and consider sequences of coefficients appearing alongside each monomial.
By lemma \ref{redlem2} we have $l(\overline{G_i}^p) \leq l(G_i)$.
Let $M$ be a monomial appearing in at least one of $\overline{G_i}^p$. 
We obtain a sequence  $(\overline{\alpha_M}^p) _{p \in S}$ of coefficients associated with $M$. Let us denote the upper bound for $A$ given in theorem \ref{stthm} by $C$, i.e.
\[ C:= B^{\sum_{i=0}^{\mu -2}D^{i}} \cdot (l(F)-1)^{\mu - 2} \prod_{j=1}^{\mu-2}\Big( l(F)^{D^j}+1\Big)^{\mu -1-i} \cdot \Big[ l(F) + (l(F)-1)\sum_{r=2}^{\mu-1} \prod_{s=1}^{r-1}(l(F)^{D^s} +1)\Big].\]

 \begin{rmk}
 For any integer $q > 2C$ we have $G = \overline{G}^q$.

\label{stabcheck}
 \end{rmk}
 \textit{Proof.} Let $M$ be an arbitrary monomial appearing in $G$ with a coefficient $\alpha_M \in \mathbb{Z}$. By theorem \ref{stthm}  we have $\overline{\alpha_M}^{q}=\alpha_M$, for every $q >2C$.
   \begin{flushright}
$\square$
\end{flushright}

 \begin{cor}
 Let $F \in \mathbb{Z}[X]^n$ be as in theorem \ref{stthm}.  Let $M_{G}$ and $M_{\overline{G}^p}$
 denote set of all monomials appearing in $G$ and $\overline{G}^p$ respectively.
 If $p > 2C$, then $M_G = M_{\overline{G}^p}$.
 \label{allmoncheck}
 \end{cor}

  By remark \ref{stabcheck} and corollary \ref{allmoncheck} if we choose $S$ in such a way that  $$\prod_{s \in S} s \geq 2C+1,$$
then using Chinese Remainder Theorem we can check that we get all monomials and that we actually observe a stabilization of all coefficients. We retrieve $G$ by  considering values obtained after stabilization  as coefficients from $\mathbb{Z}$.
The meaning of remark \ref{stabcheck} and corollary \ref{allmoncheck} is theoretical.
These observations states that the procedure can always be finished in a finite number of steps. For examples with relatively big coefficients, one can try perform reductions for some subset of prime numbers and confirm retrieving of an inverse by computing the composition of $F$ and obtained $G$.

Below we present an example which illustrates how one can use results obtained in the previous section and how this approach helps to save time and memory needed to find an inverse of a given polynomial automorphism.

\begin{ex}
 Let us consider  $F: \Q^4 \rightarrow \Q^4$ given by the following formula.
\[  F: \left\{ \begin{array}{ll}
 F_1=&X_1 \\
 F_2=&X_2 + 3 X_{1}^{3}\\
 F_3=&X_3 -159471666 X_{1}^{15} - 136514727 X_{1}^{13} - 265786110 X_{1}^{12} X_{2} - 2541294 X_{1}^{11} \\&- 182019636 X_{1}^{10} X_{2} - 177190740 X_{1}^{9} X_{2}^{2} + 972 X_{1}^{9} - 2541294 X_{1}^{8} X_{2} - 91009818 X_{1}^{7} X_{2}^{2} \\& - 
59063580 X_{1}^{6} X_{2}^{3} + 405 X_{1}^{7} + 1215 X_{1}^{6} X_{2} - 847098 X_{1}^{5} X_{2}^{2} - 20224404 X_{1}^{4} X_{2}^{3} \\ &- 9843930 X_{1}^{3} X_{2}^{4} - 93717 X_{1}^{6} X_{3} - 81 X_{1}^{6} X_{4} + 27 X_{1}^{5} + 270 X_{1}^{4} X_{2
} + 486 X_{1}^{3} X_{2}^{2} \\ &- 94122 X_{1}^{2} X_{2}^{3} - 1685367 X_{1} X_{2}^{4} - 656262 X_{2}^{5} - 27 X_{1}^{4} X_{3} - 62478 X_{1}^{3} X_{2} X_{3} \\ &- 54 X_{1}^{3} X_{2} X_{4} + 9 X_{1}^{2} X_{2} + 45 X_{1} X_{2}^{2} + 63 X_{2}^{3} - 9
 X_{1} X_{2} X_{3} - 10413 X_{2}^{2} X_{3} - 9 X_{2}^{2} X_{4} \\
 F_4=&X_4 + 184508717562 X_{1}^{15} + 158000696361 X_{1}^{13} + 307514529270 X_{1}^{12} X_{2} + 2985782067 X_{1}^{11} \\&+ 210667595148 X_{1}^{10} X_{2} + 205009686180 X_{1}^{9} X_{2}^{2} + 1691280 X_{1}^{9} + 2985782067 X_{1}^{8} X_{2} \\ & + 1053337
97574 X_{1}^{7} X_{2}^{2} + 68336562060 X_{1}^{6} X_{2}^{3} + 1216458 X_{1}^{7} + 1127763 X_{1}^{6} X_{2} \\ &+ 995260689 X_{1}^{5} X_{2}^{2} + 23407510572 X_{1}^{4} X_{2}^{3} + 11389427010 X_{1}^{3} X_{2}^{4} + 108430569 X_{1}^{6} X_{3} + \\ &93
717 X_{1}^{6} X_{4} + 810891 X_{1}^{4} X_{2} + 188082 X_{1}^{3} X_{2}^{2} + 110584521 X_{1}^{2} X_{2}^{3} + 1950625881 X_{1} X_{2}^{4} \\ &+ 759295134 X_{2}^{5} + 62478 X_{1}^{4} X_{3} + 72287046 X_{1}^{3} X_{2} X_{3} + 27 X_{1}^{4} X_{4} + 6
2478 X_{1}^{3} X_{2} X_{4} \\ & + 135135 X_{1} X_{2}^{2} + 27 X_{2}^{3} + 9 X_{1}^{2} X_{3} + 20826 X_{1} X_{2} X_{3} + 12047841 X_{2}^{2} X_{3} \\ & + 9 X_{1} X_{2} X_{4} + 10413 X_{2}^{2} X_{4}
\end{array} \right. .\]
Observe that all coefficients are integer numbers, hence there is no need to perform denominators  clearing procedure. One can perform an algorithm over $\mathbb{Z}$ and find an inverse mapping $G$.
\[  G: \left\{ \begin{array}{ll}
G_1 =& X_1 \\
G_2 =& X_2 -3 X_{1}^{3} \\
G_3 =& X_3 -9 X_{1}^{2} X_{2} - 45 X_{1} X_{2}^{2} - 63 X_{2}^{3} + 9 X_{1} X_{2} X_{3} + 10413 X_{2}^{2} X_{3} + 9 X_{2}^{2} X_{4} \\
G_4 = & X_4 -135135 X_{1} X_{2}^{2} - 27 X_{2}^{3} - 9 X_{1}^{2} X_{3} - 20826 X_{1} X_{2} X_{3} - 12047841 X_{2}^{2} X_{3} \\& - 9 X_{1} X_{2} X_{4} - 10413 X_{2}^{2} X_{4}
\end{array} \right. .\]

\begin{figure}[h]
    \centering
    \includegraphics[scale=0.42]{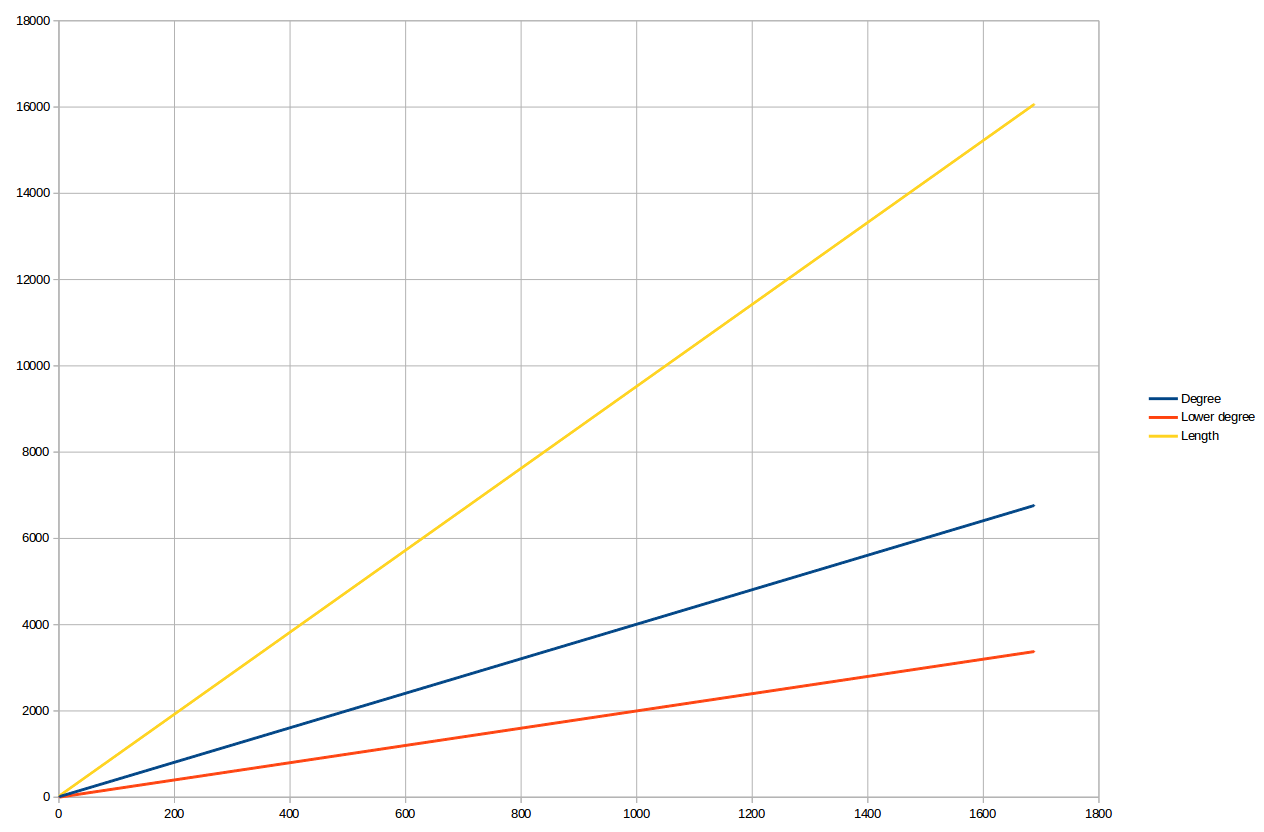}
    \caption{Degrees, lower degrees and lengths of the third coordinate of the sequence produced by algorithm for $F$.}
    \label{fig:original_plot}
\end{figure}
These calculations take 57 minutes and 32 seconds and consume 7GB RAM. According to theorem \ref{symthm}, we need to perform at most 1688 steps of the algorithm in order to find the inverse mapping. It appears that algorithm does not stop in 1688 steps for any coordinate. For the previous examples we presented degrees, lower degrees and lengths of chosen coordinate of polynomial mappings produced by the algorithm in table. In this example due to large size of the numbers we present such a data in figure \ref{fig:original_plot} instead.

Alternatively we can use reductions modulo prime numbers and then obtain $G$ using Chinese Remainder Theorem. For example, after reducing $F$ modulo 2 we obtain  $\overline{F}^2~:~\mathbb{F}_2 \rightarrow \mathbb{F}_2$ given by the following formula.
\[  \overline{F}^2: \left\{ \begin{array}{ll}
 \overline{F}^2_1=&X_1 \\
 \overline{F}^2_2=&X_2 + X_{1}^{3}\\
 \overline{F}^2_3=&X_3 + X_{1}^{13} + X_{1}^{7} + X_{1}^{6} X_{2} + X_{1}^{6} X_{3} + X_{1}^{6} X_{4} + X_{1}^{5} + X_{1} X_{2}^{4} + X_{1}^{4} X_{3} \\ &+ X_{1}^{2} X_{2} + X_{1} X_{2}^{2} + X_{2}^{3} + X_{1} X_{2} X_{3} + X_{2}^{2} X_{3} + X_{2}^{2} X_{4} \\
 \overline{F}^2_4=&X_4 + X_{1}^{13} + X_{1}^{11} + X_{1}^{8} X_{2} + X_{1}^{6} X_{2} + X_{1}^{5} X_{2}^{2} + X_{1}^{6} X_{3} + X_{1}^{6} X_{4} + X_{1}^{4} X_{2} \\ & + X_{1}^{2} X_{2}^{3} + X_{1} X_{2}^{4} + X_{1}^{4} X_{4} + X_{1} X_{2}^{2} + X_{2}^{3} + X_{1}
^{2} X_{3} + X_{2}^{2} X_{3} + X_{1} X_{2} X_{4} + X_{2}^{2} X_{4}
\end{array} \right. .\]

Algorithm allows us to find its inverse $\overline{G}^2$.
\[  \overline{G}^2: \left\{ \begin{array}{ll}
\overline{G}^2_1 =& X_1 \\
\overline{G}^2_2 =& X_2 + X_{1}^{3} \\
\overline{G}^2_3 =& X_3 + X_{1}^{2} X_{2} + X_{1} X_{2}^{2} + X_{2}^{3} + X_{1} X_{2} X_{3} + X_{2}^{2} X_{3} + X_{2}^{2} X_{4} \\
\overline{G}^2_4 = & X_4 + X_{1} X_{2}^{2} + X_{2}^{3} + X_{1}^{2} X_{3} + X_{2}^{2} X_{3} + X_{1} X_{2} X_{4} + X_{2}^{2} X_{4}
\end{array} \right. .\]
These calculations take 15 seconds using 0,57 GB RAM. The maximum number of steps given by the theorem \ref{symthm} is equal to 1099. Algorithm does not stop until then.

We observe that it is enough to consider $S= \{2,3,5,7,11,13,17,19,23,29\} $ to reconstruct $G$. One can check that by computing the composition of $F$ and obtained $G$. Times of execution and consumption of memory are presented in table \ref{tab:summary}. Degrees, lower degrees and length appearing in the third coordinate of the sequence calculated during execution of algorithm for mapping
$\overline{F}^2$ are presented in figure \ref{fig:2_plot}.
\begin{table}[!htbp]
\centering
\begin{tabular}{|c||c|c|}
    \hline
        $p$ & Time of execution & RAM used \\ \hline
        2 & 15 s & 0.57 GB \\\hline
        3 & 1 s & 0.2 GB \\\hline
        5 & 3 min 5 s & 2.96 GB \\\hline
        7 & 5 min & 3.94 GB \\\hline
        11 & 9 min 21 s & 5.58 GB \\\hline
        13 & 4 min 13 s &  1.78 GB\\\hline
        17 & 9 min 40 s & 5.65 GB \\\hline
        19 & 10 min 32 s & 5.63 GB \\\hline
        23 & 11 min 6 s & 6.1 GB \\\hline
        29 &    12 min 59 s & 6.53 GB \\\hline
    \end{tabular}
    \caption{Execution time and memory consumption.}
    \label{tab:summary}
\end{table}

One can observe that the reduction approach run in sequence can last longer than the direct approach for the mapping $F \in \mathbb{Q}[X]^4$. However the memory consumption can be smaller. Hence this method can be used for computer with less amount of memory installed. 
This observation allows us also to use reduction approach together with parallel computations. In this way the time of inverting $F$ can be significantly reduced.

\begin{figure}
    \centering
    \includegraphics[scale=0.48]{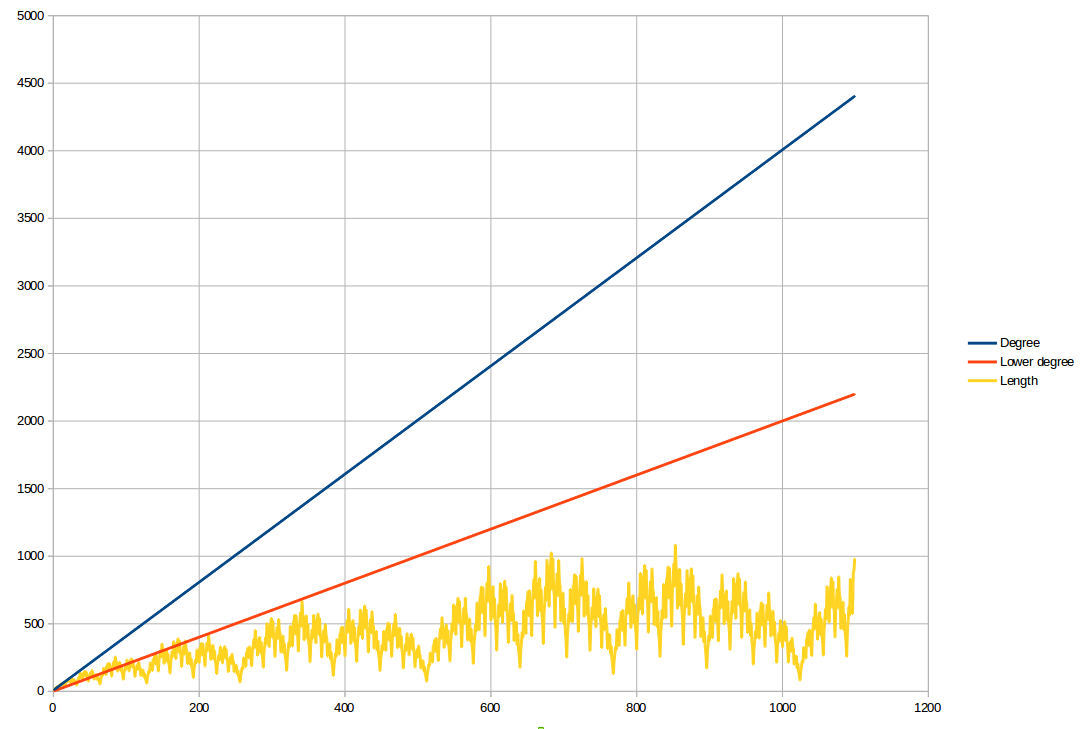}
    \caption{Degrees, lower degrees and lengths of the third coordinate of the sequence produced by algorithm for $\overline{F}^{2}$.}
    \label{fig:2_plot}
\end{figure}
\end{ex}

\end{document}